\documentclass[greek,english,12pt,reqno]{amsart}
\usepackage{babel}
\usepackage{amssymb}
\usepackage{amsxtra}
\usepackage[mathscr]{eucal}
\usepackage{graphics}
\usepackage{epsfig}

\newtheorem{theorem}{Theorem}[section]

\newtheorem{proposition}[theorem]{Proposition}

\theoremstyle{definition}

\newtheorem{definition}[theorem]{Definition}
\newtheorem{conjecture}[theorem]{Conjecture}

\def\cR{{\mathcal R}}
\def\cS{{\mathcal S}}
\def\cT{{\mathcal T}}
\def\cU{{\mathcal U}}
\def\ca {\mathcal{A}}
\def\beq {\begin{equation}}
\def\endq {\end{equation}}
\newcommand{\twolinesum}[2]{\sum_{\substack{{\scriptstyle #1}\\
{\scriptstyle #2}}}}

\newcommand{\gmod}[1]{\, ({\rm{mod}} \, #1)}

\providecommand{\Real}{\mathop{\rm Re}\nolimits}%
\providecommand{\Imag}{\mathop{\rm Im}\nolimits}%

\begin{document}
\title{On primitive divisors of $n^2+b$}
\author{Graham Everest and Glyn Harman}
\address{(GE) School of Mathematics, University of East Anglia,
Norwich NR4 7TJ, UK} \address{(GH) Department of Mathematics,
Royal Holloway, University of London, Egham, Surrey TW20 0EX UK}

\email{g.everest@uea.ac.uk} \email{g.harman@rhul.ac.uk}

\begin{abstract}
We study primitive divisors of terms of the sequence $P_n=n^2+b$,
for a fixed integer~$b$ which is not a negative square. It seems
likely that the number of terms with a primitive divisor has a
natural density. This seems to be a difficult problem. We survey
some results about divisors of this sequence as well as provide
upper and lower growth estimates for the number of terms which
have a primitive divisor.
\end{abstract}

\maketitle

\section{Primitive prime divisors} Given $b$, an integer which is not a
negative square, consider the integer sequence with $n$th term
$P_n=n^2+b.$ It seems likely \cite{Bat} that infinitely many of
the terms are prime but a proof seems elusive. Perhaps this
mirrors the status of the Mersenne Prime Conjecture, which
predicts that the sequence with $n$th term $M_n=2^n-1$ contains
infinitely many prime terms. At least with the Mersenne sequence,
an old result shows that primes are produced in a less restrictive
sense.
\begin{definition}
Let $(A_n)$ denote a sequence with integer terms. We say an
integer $d>1$ is a {\it primitive divisor} of $A_n$ if
\begin{enumerate}\item $d|A_n$ and \item gcd$(d,A_m)=1$ for all non-zero
terms~$A_m$ with $m<n$.\end{enumerate}\end{definition}

In 1886 Bang \cite{bang} showed that if~$a$ is any fixed integer
with $a>1$ then the sequence with $n$th term~$a^n-1$ has a
primitive divisor for any index~$n>6$. This is remarkable because
the number 6 is uniform across all $a$ and it is small. Before we
say any more about polynomials, a short survey follows indicating
the incredible influence of Bang's Theorem.

\subsection{Primitive divisor theorems}
In 1892 Zsigmondy obtained the generalization that for any choice
of $a$ and $b$ with $a>b>0$, the term~$a^n-b^n$ has a primitive
divisor for any index~$n>6$. This lovely result was re-discovered
several times in the early 20th century and it has turned out to
be quite applicable. See \cite{cp} and the references therein
where applications to Group Theory are discussed. For example, the
order of the group $GL_n(\mathbb F_q)$ has a primitive divisor for
all large~$n$. Thus Sylow's Theorem can be invoked to deduce
information about the structure of the group.

The next major theoretical advance was made by Carmichael. Let~$u$
and~$v$ denote conjugate quadratic integers; in other words, zeros
of a monic irreducible polynomial with integer coefficients.
Consider the integer Lucas sequence defined by
$$
U_n=(u^n-v^n)/(u-v).
$$
The Fibonacci sequence $(F_n)$ arises from the roots of the
polynomial $x^2-x-1$. Carmichael \cite{carmichael} showed that if
$u$ and $v$ are real then $U_n$ has a primitive divisor for
$n>12$. This is a sharp result because $F_{12}$ does not have a
primitive divisor. Less is currently known about the corresponding
Lehmer-Pierce sequence
$$
V_n=(u^n-1)(v^n-1).
$$
Kalman Gy\"ory pointed out to the first author that if $uv=1$ then
$V_n$ has a primitive divisor for all~$n$ beyond some (actually
uniform) bound; on the other hand, if $uv=-1$ then $V_{2k}$ does
not have a primitive divisor if $k$ is odd, because
$V_{2k}=-V_k^2$. This second observation is actually quite germane
to this paper; see Theorem \ref{maincorollary}. In fact the set of
terms with a primitive divisor has natural density equal to
$\frac34$ (cf. Conjecture \ref{conjecture}). At the conference,
Richard Pinch remarked that certain Lehmer-Pierce sequences count
orders of groups: this time the groups are $E(\mathbb F_{p^n})$,
where $E$ denotes an elliptic curve.

Bilu, Hanrot and Voutier~\cite{MR2002j:11027} used powerful
methods from Diophantine analysis to prove, in the general case,
that~$U_n$ has a primitive divisor for any~$n>30$. Again this is a
sharp result as the sequence generated by the polynomial $x^2-x+2$
illustrates. Finally, Silverman \cite{silabc} obtained a primitive
divisor theorem for Elliptic Divisibility Sequences and a uniform
version appears in \cite{emw} for a certain class of sequences.

\subsection{Primitive divisors of $n^2+b$}

\begin{theorem}\label{maincorollary}
Infinitely many terms of the sequence~$n^2+b$ do not have a
primitive divisor.
\end{theorem}

The proof of Theorem \ref{maincorollary} follows very easily from
a result of Schinzel \cite{schinzel} and will be discussed
shortly. Schinzel's proof manufactures a very thin set of terms
with no primitive divisor. Dartyge \cite{Dartyge} has improved
Schinzel's result for $n^2+1$ (and in principle the method works
for $n^2 + b$ also). The aim of this paper to obtain a better
grasp on the set of terms with no primitive divisor. We will also
consider whether the set of indices $n$ for which~$P_n$ has a
primitive divisor has a natural density. Apparently this lies
quite deep.

For other interesting approaches to the study of divisors of
quadratic integral polynomials; consult \cite{Dartyge},
\cite{Des}, \cite{MR1324141}, \cite{MR0163874}, \cite{MR0204383},
\cite{mckee1}, \cite{mckee2} and \cite{MR1776618}. For higher
order polynomials there is also the paper \cite{DartMart}.

\subsection{The greatest prime factor}

Let $P^+(m)$ denote the greatest prime factor of the integer
$m>1$. There is a wealth of literature about $P^+(n^2+b)$
concerned with the fact that $P^+(n^2+b)\rightarrow \infty$ as
$n\rightarrow \infty$, see \cite[Chapter 7]{st}. In a slightly
different direction, Luca \cite{luca} has recently revived an old
method of Lehmer's \cite{lonstormer} to show that, given~$B$, the
set of indices for which $P^+(n^2+1)<B$ is efficiently computable.
Carmichael's result mentioned earlier for Lucas sequences plays a
key role. He illustrates his method by showing that when $B=101$,
$n \le 24208144$.

The following is an easy proposition, see \cite{candt} or
\cite{AMM}, which relates $P^+(n^2+b)$ to the existence of a
primitive divisor.

\begin{proposition}\label{primeishiff} For all~$n>|b|$,
the term~$P_n=n^2+b$ has a primitive divisor if and only
if~$P^+(n^2+b)>2n$. For all~$n>|b|$, if~$P_n$ has a primitive
divisor then that primitive divisor is a prime and it is unique.
\end{proposition}

\begin{proof}[Proof of Theorem~\ref{maincorollary}]
Results of Schinzel~\cite[Th.~13]{schinzel} show that for
any~$\alpha>0$, $P^+(n^2+b)$ is bounded above by~$n^{\alpha}$ for
infinitely many~$n$. Taking $\alpha=\frac{1}{2}$,
Proposition~\ref{primeishiff} shows that~$P_n=n^2+b$ fails to have
a primitive divisor infinitely often.
\end{proof}

Given $x>1$, Schinzel's method constructs fewer than $\log x$
terms $P_n$ with $n<x$ having no primitive divisor. For $\alpha
> \frac{149}{179}$, Dartyge \cite{Dartyge} showed that
\[
\left|\{n \le x: P^+(n^2+1) < x^{\alpha} \}\right| \gg x.
\]
It should be noted that the implied constant is very small,
involving, as it does, a term $2^{-\delta^{-2}}$ where $\delta$
``est extr\^ement petit" \cite[p.3 line 10]{Dartyge}. In this
paper we prove the following, which provides good upper and lower
estimates for the number of terms with a primitive divisor.

\begin{theorem}\label{bb} Supposing~$-b$ is not an integer square, define
$$
\rho_b(x)=\left|\{n\le x:n^2+b \mbox{ has a primitive divisor
}\}\right|.
$$
For all sufficiently large $x$ we have
\[
0.5324 < \frac{\rho_b(x)}{x} < 0.905.
\]
\end{theorem}

\subsection{Natural density}
Integers $m$ with the property $P^+(m)>2\sqrt m$ were studied by
Chowla and Todd \cite{candt}. They proved that the set of these
numbers has natural density $\log 2$. Perhaps this suggests the
following:

\begin{conjecture}\label{conjecture} If~$-b$ is not an integer
square then $\rho_b(x)\sim x\log 2.$
\end{conjecture}

With the availability and power of modern computers, one would
usually resort to some computational evidence in support of such a
conjecture. The authors of \cite{AMM} looked for such evidence.
Whilst they found nothing to clearly contradict the conjecture,
neither did they find overwhelming evidence to support it. The
problem is that the convergence to the natural density is very
slow.

The reason for this might best be explained as follows. Chowla and
Todd's proof uses Mertens' Theorem about the asymptotic formula
for the sum of inverse primes:
$$\sum_{p<x}\frac{1}{p}=\log \log x + C + O\left(\frac{1}{\log x}\right).$$
The main term of this formula grows very slowly and the error term
shrinks very slowly as well. Perhaps, somehow, this lies behind
the extremely slow convergence to the natural density of terms
with primitive divisor, as in Conjecture \ref{conjecture}. In
addition, the arithmetical nature of the sequence $n^2+b$ plays a
significant r\^ole when discussing its very large prime divisors
(see \eqref{polysieve} below) and this will affect what happens
for `small' $x$. Our paper concludes with an explanation as to why
we are not holding our breath about a proof of Conjecture
\ref{conjecture}.

\section{Simple bounds}

The article \cite{AMM} gives some simple estimates for
$\rho_{b}(x)$ which are sketched below. These are recalled here as
a way in to the harder methods. The first bound in
(\ref{firsttheorembound}) counts indices which produce no
primitive divisor. It is much better than the bound obtained from
\cite{schinzel} but the set of indices still has density zero and
perhaps indicates the limit of elementary methods. The second
bound in (\ref{firsttheorembound}) is very easy but already gives
a good estimate for the density of terms with a primitive divisor
if it exists.

\begin{theorem}\label{infinitelymany}
For all sufficiently large $x$,
\begin{equation}\label{firsttheorembound} \frac{x}{\log x} \ll
x-\rho_{b}(x)\mbox{ and } \tfrac 12 x - \rho_b(x)\ll \frac{x}{\log
x}.
\end{equation}
\end{theorem}

The proofs use little apart from well-known estimates for sums
over primes, which can be found in the book of
Apostol~\cite{MR0434929}.  Both begin with an old idea of
Chebychev which is used frequently as the starting point of
investigating the greatest prime factor of certain sequences (see
\cite[Chapter 2]{Hooley} for example).

Apart from a finite number of primes, any prime~$p$ that
divides~$n^2+b$ has the property that~$-b$ is a quadratic residue
modulo~$p$. Let~$\cR$ denote the set of odd primes for which~$-b$
is a quadratic residue; notice that~$\cR$ comprises the
intersection of a finite union of arithmetic progressions with the
set of primes. Write
$$
Q_x=\prod_{n=1}^x\vert P_n\vert
$$
and denote by~$\omega (Q_x)$ the number of prime divisors
of~$Q_x$. By Proposition~\ref{primeishiff} it is sufficient to
bound~$\omega (Q_x)$ because, with finitely many exceptions, a
primitive divisor is unique.

\subsection{Sketch proof of Theorem \ref{infinitelymany}}

Define
$$
\cS=\{p\in \cR: p|Q_x, p<2x\} \mbox{ and } \cS'=\{p\in \cR: p|Q_x,
p \ge 2x\}.
$$
Let~$s=\vert \cS\vert$ and~$s'=\vert \cS'\vert$. We seek bounds
for~$s+s'$. By Dirichlet's Theorem on primes in arithmetic
progression it is sufficient to estimate~$s'$. Following
Chebychev's method, use Stirling's Formula to obtain
\begin{equation}\label{estimate1}
\sum_{p|Q_x}e_p\log p = \log Q_x  = 2x\log x +O(x)
\end{equation}
where the left-hand side corresponds to the prime decomposition
of~$Q_x$, for positive integers~$e_p$. The sum on the left-hand
side of \eqref{estimate1} decomposes according to the definitions
of~$\cS$ and~$\cS'$ to give
\begin{equation}\label{estimate2}
\sum_{p\in \cS}e_p\log p +\sum_{p\in \cS'}\log p=\log Q_x,
\end{equation}
noting that~$e_p=1$ whenever~$p\ge 2x$. It is easy to show that
\begin{equation}\label{sumoverpinS}\sum_{p\in \cS}e_p\log p = x\log
x + O(x).
\end{equation} Combining~\eqref{estimate1}, \eqref{estimate2}
and~\eqref{sumoverpinS} gives
\begin{equation}\label{protocheb}
x\log x + O(x) = \sum_{p\in \cS'}\log p.
\end{equation}
The right hand side is bounded above by $s'\log (x^2+1)$ yielding
a lower bound for~$s'$.

The second bound in (\ref{firsttheorembound}) arises similarly
using a finer partition of the set~$\cS'$
\begin{eqnarray*}
\cT&=&\{p\in \cR: p|Q_x, 2x<p<Kx\};\\
\cU&=&\{p\in \cR\mid p|Q_x, Kx<p\},
\end{eqnarray*}
for~$K>2$. Write~$t=|\cT|$ and~$u=|\cU|$ then we seek an upper
bound for expression~$t+u$. Using the definitions of~$\cT$
and~$\cU$ as well as equation~\eqref{sumoverpinS} shows that
$$
\sum_{p\in \cT}\log p +\sum_{p\in \cU}\log p = x\log x + O(x).
$$
The extra leverage comes because the left-hand side is greater
than
$$
t\log x + u\log (Kx).
$$
Now $K$ can be chosen judiciously to beat the other $O$-constants.
An upper bound for $t+u$ follows easily and hence the second bound
in (\ref{firsttheorembound}).

\medskip
\noindent {\bf Note} Actually $K$ can be taken as large as $\log
x$ which yields
$$\frac{x\log \log x}{\log x} <
x-\rho_{b}(x)
$$
for all large $x$. But this still fails to produce a positive
density set.

\section{Better bounding}
It is the aim of this section to prove Theorem \ref{bb}. Take $b =
1$ for simplicity, so we can drop the subscript $b$ on $\rho$; as
with \cite{Des} the arguments in \cite{Hooley} can be used to
generalise to $b \ne 1$.  We then prove the following.
\begin{theorem}
For all sufficiently large $x$ we have
\[
0.5324 < \frac{\rho(x)}{x} < 0.905.
\]
\end{theorem}

Let
\[
N_x(p) = \twolinesum{x \le n < 2x}{p|n^2 + 1} 1.
\]
The previous section shows it is sufficient to estimate
$$\sum_{p \ge 2x}N_x(p).
$$
Re-casting (\ref{protocheb}) using this definition:
\begin{equation}\label{final}
\sum_{p \ge 2x} N_x(p) \log p = x \log x\ + O(x).
\end{equation}
The extreme cases arise if most of the contribution to this sum
comes from $p$ around $2x$ in size, or around $4x^2$ in size.  In
the former case the bound $\log p \ge \log x$ gives the trivial
bound \beq\label{upper} \sum_{p \ge 2x} N_x(p) < x,
\endq
which is weaker than the first bound of the last section. On the
other hand, $\log p \le 2 \log x + O(1)$ gives \beq\label{lower}
\sum_{p \ge 2x} N_x(p) > \tfrac12 x + o(x),
\endq
which is essentially the second bound of the last section.

We could obtain improved results if we had better information about
the following expression:
\[ V_x(v) = \sum_{v < p \le ev} N_x(p).
\]
It is a good exercise to show that \beq\label{hope} V_x(v) \sim
\frac{x}{\log v}
\endq
implies the conjecture.  Unfortunately, the asymptotic formula
\eqref{hope} is not expected to be true for very large $v$, in
view of the arithmetic nature of $n^2 + 1$ (see below). However,
it is expected that \eqref{hope} will be true for $v < x^{2 -
\epsilon}$ for any $\epsilon > 0$ and this suffices to prove the
conjecture.

\subsection{A better upper bound for $\rho(x)$}

We begin by modifying the definitions to allow us to use the
Deshouillers-Iwaniec method in \cite{Des}. To be precise we must
use smooth functions in order to apply the mean-value estimates in
\cite{Des} for Kloostermann sums. Let $\epsilon, \eta$ be two
small positive quantities. Let $b(u)$ be a function satisfying
$b(u) \in [0,1]$ for all $u \in \mathbb{R}$, with
\[
\begin{split}
b(u) &= \begin{cases} 1  &\text{if} \ (1 + \epsilon)x \le u \le (2-\epsilon)x\\
0  &\text{if}  \ u \le x \  \text{or} \ u \ge 2x,
\end{cases}\\
\frac{{\rm d}^r b(u)}{{\rm d}u^r}  &\ll_{r,\epsilon} u^{-r} \ \text{for all} \ r \in \mathbb{N}.
\end{split}
\]
We redefine $N_x(p)$ to be
\[
N_x(p) = \twolinesum{x \le n < 2x}{p|n^2 + 1} b(n).
\]
An upper bound for this summed over $p$ will give us an upper bound for the original problem,
since the two quantities will differ by at most
\[
\tfrac32 \epsilon x.
\]
Now write
\[
X = \int_{x}^{2x} b(u) \, du, \quad \left| \ca_d \right| = \sum_{n^2 + 1 \equiv 0 \gmod{d}} b(n).
\]
By the working on \cite[p.2]{Des} we can modify the Chebychev argument to give
\[
\sum_{p} \left|\ca_p \right| \log p = 2X \log x + O(x).
\]
Also, as shown in \cite{Des}, we have
\[
\sum_{p \le x} \left| \ca_p\right| = X \log x + O(x).
\]
Let
\[
P_x = \max_{|\ca_p| \ne 0} p \ = \ x^{\sigma} \ \ \text{say.}
\]
We therefore have
\[
\sum_{x \le p \le P_x} \left| \ca_p \right| \log p = X \log x + O(x).
\]
Deshouillers and Iwaniec then estimate this sum as
\[
\sum_{1 \le j \le J} S(X,V_j) + O(x),
\]
where $V_j = 2^j x$ and
\[
S(x,V_j) = \sum_{V_j < p \le 4V_j} C_j(p) \log p.
\]
Here the infinitely differentiable functions $C_j(u) \in [0,1]$ are supported in $[V_j, 4V_j)$, with
\[
\sum_{1 \le j \le J} C_j(u) = \begin{cases} 1 &\text{if} \ 2x < u \le P_x\\
0 &\text{if} \ u < x \ \text{or} \ u > P_x.
\end{cases}
\]
After several transformations and an application of the Rosser-Iwaniec sieve in tandem with their
own sophisticated mean-value estimate for averages of Kloostermann sums, they prove that
\beq\label{DIineq}
S(x,V_j) \le \frac{2}{\log D_j} \int C_j(u) \frac{\log u}{u} \, du \left(1 + O\left( \frac{1}{\log D} \right)  \right).
\endq
Here $D_j = x^{1 - \eta} V_j^{-\frac12}$.  From this they deduce that $\sigma$ is not less than the solution to
\[
2 - \sigma - 2 \log(2 - \sigma)  = \tfrac54.
\]
That is, $\sigma = 1.202468 \ldots$.

Now, the worst case scenario for the upper bound \eqref{upper} is
if \eqref{DIineq} holds with equality for each $V_j$.  This gives
\[
\begin{split}
\sum_{p \ge 2x} N_x(p) &\le
\sum_{1 \le j \le J} S(X,V_j)(\log V_j)^{-1} + O(x (\log x)^{-1})\\
&=x \left(1 + O(\log x)^{-1}\right) \int_1^{\sigma} \frac{2}{1 - t/2} \, dt\\
 &= (2 \sigma - \tfrac32)x\left(1 + O((\log x)^{-1})\right)   < 0.905 x.
\end{split}
\]

\subsection{A better lower bound for $\rho(x)$} Now we need to show
that not all the contribution comes from primes near $x^2$.  This
is a relatively simple application of an upper bound sieve to the
set
\[
\{m: m \ell = n^2+b, x < n \le 2x\} \ \ \text{for $\ell$ in some range}.
\]
In this case we can apply the sieve with {\it distribution level}
\[
D_{\ell} = \frac{x}{\ell (\log x)^A}
\]
for some $A$, by an elementary argument: this corresponds
to $V_j/x$ in the last section.  Of course, this is why the elementary argument
is no good for $V_j$ near $x$ in size.
The crossover point between the two methods is at $V_j = x^{\frac43}$, but we can get
nowhere near this value
for the problem discussed in \cite{Des}.
For $\ell = 1$ the problem is the well-known one of representing almost-primes by values
of $n^2 + 1$ and giving an
upper bound for the number of prime values of this polynomial.  By \cite[Theorem 5.3]{Hal}
(or see \cite[p.66]{Greaves})we have
\beq\label{polysieve}
\twolinesum{x \le n \le 2x}{n^2 + 1 = p} 1 \le
\frac{2x}{\log x} \prod_p \left(1 - \frac{\chi(p)}{p-1}  \right)\left(1 +
O\left(\frac{\log \log 3x}{\log x}  \right)  \right).
\endq
Here $\chi(n)$ is the non-trivial character $\gmod{4}$. Note the important
product over primes above which encodes
 arithmetical information relevant to the polynomial $n^2 + 1$.  This did not
 arise in the previous section since summing over a sufficiently
long range for $\ell$ smooths out this factor (compare \cite[\S 8]{Des}).
It is expected that \eqref{polysieve} holds with equality if the factor $2$
is replaced by $\frac12$ on the right-hand side, see \cite{Hardy,Bat}.

We obtain our desired bound by first considering
\[
W(L,x) = \sum_{L \le \ell \le eL} \twolinesum{x \le n \le 2x}{n^2 + 1 = \ell p} 1.
\]
Write $\omega(d)$ for the number of solutions to $n^2 + 1 \equiv 0
\gmod{d}$ and let $\{\lambda_d\}_{d \le D}$ be the Rosser upper
bound sieve of level $D_L = x(L (\log x)^A)^{-1}$ as described in
\cite[\S 4]{Des} and explicitly constructed in \cite[Chapter
4]{Greaves}.  We then have
\[
\begin{split}
W(L,x) &\le
\sum_{L \le \ell \le eL} \sum_{d \le D_L} \lambda_d \twolinesum{x \le n \le 2x}{n^2 + 1 \equiv 0 \gmod{d\ell}} 1\\
&=
\sum_{L \le \ell \le eL} \sum_{d \le D_L} \lambda_d \omega(d \ell) \left(\frac{x}{d \ell}
+ O(1)  \right)\\
&= x \sum_{L \le \ell \le eL} \sum_{d \le D_L} \lambda_d
\frac{\omega(d \ell)}{d \ell} + O\left(LD_L (\log x)^2  \right).
\end{split}
\]
In the above we have noted that $\omega(d \ell) \le \tau(d)
\tau(\ell)$ and used the well-known average value of the divisor
function $\tau(n)$ to give a bound for the error term. We then use
a similar analysis to that in \cite[\S 8]{Des} to produce the
`main term'.  We give all the details that differ from \cite{Des}
here for completeness.

Firstly write
\[
\sum_{L \le \ell \le eL} \sum_{d \le D_L} \lambda_d \frac{\omega(d \ell)}{d \ell}
= \sum_{d \le D_L} \lambda_d \frac{\omega(d)}{d} J(d,L).
\]
Now put
\[
L(s,d) = \sum_{m=1}^{\infty} \frac{\omega(dm)}{\omega(d)m^s}.
\]
Note by \cite[Lemma 4]{Des} that
\[
L(s,d) = \frac{\zeta(s) L(s,\chi)}{\zeta(2s)} \prod_{p|d}\left(1 + \frac{1}{p^s} \right)^{-1}.
\]
Using Perron's formula (\cite[Theorem 3.12]{Titch}) with $T=x, c =
(\log x)^{-1}$ we have
\[
J(d,L) = \frac{1}{2 \pi i} \int_{c-ix}^{c+ix} L(s+1,d) \frac{(eL)^s - L^s}{s} \, ds +
O\left(x^{-\frac12}  \right).
\]
The final term is negligible. (Actually, this has been estimated
very crudely, in reality it is $O(x^{\epsilon -1})$). Now take the
contour of integration back to $\Real s = - \frac12$. The pole at
$s=0$ gives a term
\[
\frac{L(1,\chi)}{\zeta(2)}
\prod_{p|d}\left(1 + \frac{1}{p} \right)^{-1}.
\]
The pair of integrals on $\Imag s = \pm x$ give a
negligible contribution ($O(x^{\epsilon -\frac23})$), using
\[
\max(|\zeta(s)|, |L(s,\chi)|) \ll T^{\frac16} \ \ \text{for} \ 1 \le |\Imag s | \le T, \Real s \ge \tfrac12.
\]
The integral on
the new contour can be estimated using:
\[
\int_{-T}^T \frac{|\zeta(\tfrac12 + it)|^2}{1 + |t|} \, dt \ll (\log T)^2,
\]
with the same bound applying when $\zeta(s)$ is replaced by $L(s,\chi)$,
together with (\cite[p.135]{Titch})
\[
\frac{1}{\zeta(1 + it)} \ll \log T \ \ \text{for} \ |t| \le T,
\]
and
\[
\left|\prod_{p|d}\left(1 + \frac{1}{p^{s+1}} \right)^{-1}\right| \le
\prod_{p|d}\left(1 - \frac{1}{p^{\frac12}} \right)^{-1} < \tau(d).
\]
This gives a bound for the integral which is
\[
\ll \frac{\tau(d)(\log x)^3}{L^{\frac12}}.
\]

Thus
\[
\sum_{L \le \ell \le eL} \sum_{d \le D_L} \lambda_d \frac{\omega(d \ell)}{d \ell}
= \sum_{d \le D_L} \lambda_d \frac{\omega'(d)}{d} + O\left(\sum_{d \le D_L}
\frac{\omega(d)\tau(d)(\log x)^3}{L^{\frac12} d} \right).
\]
Here
\[
\omega'(d) = \omega(d) \prod_{p|d}\left(1 + \frac{1}{p} \right)^{-1}.
\]
The rest of the working follows {\it mutatis mutandis} from
\cite[p.10]{Des}. Hence
\[
x \sum_{L \le \ell \le eL} \sum_{d \le D_L} \lambda_d \frac{\omega(d \ell)}{d \ell}
= \frac{2x}{\log D_L} \left(1 + O\left(\frac{1}{\log D_L}  \right)  \right) +
O\left(\frac{x (\log x)^7}{L^{\frac12}} \right).
\]
The reader can thus see that the extra error term $O(x (\log x)^7
L^{-\frac12})$ (the log power could be reduced here by more
careful working) corresponds to the averaging over $\ell$
smoothing out the influence of the product in \eqref{polysieve},
and this must dominate the main term for small $L$ since the `main
term' will be incorrect in this case. Assuming that $D_L = xL^{-1}
(\log x)^{-4}$ and $x^{\frac34} >L > (\log x)^{18}$ we obtain
\beq\label{elupper} W(N,L) \le \frac{2x}{\log D_L} +
O\left(\frac{x}{(\log x)^2} \right).
\endq

For $L \le (\log x)^{18}$ we can establish a slightly cruder upper
bound as follows. For each value of $\ell$ we do not sieve by
primes dividing $\ell$.  This makes the $\lambda_d$ depend on
$\ell$, but we have $\lambda_d = 0$ if $(d,\ell)>1$.  Hence we can
write $\omega(d \ell) = \omega(d) \omega(\ell)$. Following the
analysis above, the remainder term remains $O(L D_L (\log x)^2)$.
The `main term' for the upper bound is now
\[
\frac{2x}{\log D_L} \sum_{L \le \ell \le eL} \frac{\omega(\ell)}{\phi(\ell)}\prod_{p \nmid \ell}
\left(1 - \frac{\chi(p)}{p-1}  \right)
\le \frac{K x}{\log D_L}
\]
for some absolute constant $K$.  The contribution from the terms with $L \le (\log x)^{18}$
 is thus $\ll (\log \log x)(\log x)^{-1}$ times the
total contribution for larger $L$.  These terms may therefore be neglected asymptotically.

Now the worst case scenario for \eqref{upper} has equality in
\eqref{elupper} for
\[x^{2-\theta} \ge L \ge (\log x)^{18}, \mbox{
where } \int_{\theta}^2 \frac{2t}{t - 1} \, dt = 1.
\]
In other words, $\theta$ is the solution to
\[
2(2 - \theta) - 2 \log(\theta - 1) = 1.
\]
This is the limit of the sequence
\[
a_1 = 2, \quad a_{n+1} = \frac12 \left(\frac32 + a_n - \log(a_n -1) \right) \ (n \ge 1),
\]
quickly giving the value $1.766249 \ldots$.
We then calculate
\[
\int_{\theta}^2 \frac{2}{t - 1} \, dt = 2 \theta - 3 > 0.5324
\ldots
\]

\section{Some implications of Conjecture \ref{conjecture}}
The following argument shows that we do not expect Conjecture
\ref{conjecture} to be settled in the near future. In the previous
section we have used the tools that have been developed for the
investigation of the greatest prime factor of $n^2 + 1$ to obtain
(rather weak) approximations to the conjecture. Now we assume the
conjecture and demonstrate that it would lead to a phenomenal
improvement for the greatest prime factor problem.

The conjecture leads to
\[
\sum_{p \ge 2x} N_x(p) \sim x \log 2.
\]
By the Chebychev argument (\ref{final}), on average in these sums,
\[
\frac{\log p}{\log x} \sim \frac{1}{\log 2} = \sigma \
\text{(say)} \ = 1.4416\ldots
\]
Hence the greatest prime factor of $n^2 + 1$ infinitely often
exceeds $n^{\sigma}$. This more than doubles the improvement of
Deshouillers-Iwaniec over the trivial estimate! However, we can do
still better using the elementary bound from the last section. The
worst case scenario now has all the contribution to the left hand
side of \eqref{final} coming from $p$ close to $x^{\sigma}$. Since
the bounds of the last section must hold (and they are better than
the Deshouillers-Iwaniec estimates in this region), this
corresponds to finding $\alpha < \sigma < \beta$ with
\[
\int_{\alpha}^{\beta} \frac{2}{t - 1} \, dt = \log 2, \quad
\int_{\alpha}^{\beta} \frac{2t}{t - 1} \, dt = 1.
\]
A little bit of manipulation gives the solution to be
\beq\label{b} \beta= 1 + \frac{1 - \log 2}{2 - \sqrt{2}}  =
1.52383 \ldots
\endq

This gives the following result.
\begin{theorem}
If Conjecture \ref{conjecture} is true, then infinitely often the
greatest prime factor of $n^2 + 1$ exceeds $n^{\beta}$ where
$\beta$ is given by \eqref{b}.
\end{theorem}


\end{document}